\documentclass[12pt,a4paper]{article}
\usepackage{amssymb}
\usepackage[latin1]{inputenc}
\usepackage{amsthm}
\usepackage{amsfonts}
\usepackage{mathrsfs}
\usepackage{graphicx}
\usepackage{amsmath}

\newtheorem{theorem}{Theorem}

\newtheorem{conjecture}[theorem]{Conjecture}
\newtheorem{corollary}[theorem]{Corollary}

\newtheorem{lemma}[theorem]{Lemma}

\newtheorem{proposition}[theorem]{Proposition}
\newtheorem{remark}[theorem]{Remark}

\begin{document}
\title{Weak analytic hyperbolicity of complements of generic surfaces of high
degree in projective 3-space\thanks{
Complements of projective hypersurfaces; Kobayashi hyperbolicity;
entire curves.}}

\author{Erwan Rousseau}
\date{}
\maketitle

\begin{abstract}
In this article we prove that every entire curve in the complement of a
generic hypersurface of degree $d\geq 586$ in $\mathbb{P}_{\mathbb{C}}^{3}$
is algebraically degenerated i.e there exists a proper subvariety which
contains the entire curve.
\end{abstract}

\section{Introduction}

A complex manifold $X$ is hyperbolic in the sense of S. Kobayashi if the
hyperbolic pseudodistance defined on $X$ is a distance (see, for example,
\cite{Ko98}). The hyperbolicity problem in complex geometry studies the
conditions for a given complex manifold $X$ to be hyperbolic. In the case of
hypersurfaces in $\mathbb{P}^{n}$ we have the Kobayashi conjectures \cite
{Ko70}:

\begin{conjecture}
\label{con1}A generic hypersurface $X\subset \mathbb{P}^{n+1}$ $(n\geq 2)$
of degree $\deg X\geq 2n+1$ is hyperbolic.
\end{conjecture}

\begin{conjecture}
\label{c2}$\mathbb{P}^{n}\backslash X$ $(n\geq 2)$ is hyperbolic for a
generic hypersurface $X\subset \mathbb{P}^{n}$ of degree $\deg X\geq 2n+1.$
\end{conjecture}

A new approach which could lead to a positive result for
conjecture \ref{con1} has been described by Y.-T. Siu in \cite
{SY04} for a bound $\delta _{n}\gg n$ on the degree. If we are
interested in the lower bound on the degree, conjecture \ref{con1}
is recently proved in \cite{Paun05} for $n=2$, $d\geq 18$ and in
\cite{Rou06} we proved a weak
form of conjecture \ref{con1} for $n=3$:

\bigskip

\textbf{Theorem (\cite{Rou06}). }\textit{For $X\subset \mathbb{P}_{\mathbb{C}%
}^{4}$ a generic hypersurface such that $d=deg(X)\geq 593,$ every
entire curve $f:\mathbb{C}\rightarrow X$ is algebraically
degenerate, i.e there exists a proper subvariety $Y\subset X$ such
that $f(\mathbb{C)}\subset Y.$}

\bigskip

Here we study the logarithmic conjecture \ref{c2} (proved for $n=2$ and $%
d\geq 15$ in \cite{E.G}) and prove the following result, which is a weak
form of the conjecture for $n=3$:

\begin{theorem}
\label{t2}For $X\subset \mathbb{P}_{\mathbb{C}}^{3}$ a generic hypersurface
such that $d=deg(X)\geq 586,$ every entire curve $f:\mathbb{C}\rightarrow
\mathbb{P}_{\mathbb{C}}^{3}\backslash X$ is algebraically degenerated i.e
there exists a proper subvariety $Y\subset \mathbb{P}_{\mathbb{C}}^{3}$ such
that $f(\mathbb{C)}\subset Y.$
\end{theorem}

The proof is based on two techniques.

The first one is a generalization in the logarithmic setting of an approach
initiated by Clemens \cite{Cle}, Ein \cite{Ein}, Voisin \cite{Voi} and used
by Y.-T. Siu \cite{SY04} to construct vector fields on the total space of
hypersurfaces in the projective space. Here we construct vector fields on
logarithmic spaces.

The second one is based on bundles of logarithmic jet differentials (see
\cite{DL96}). The idea, in hyperbolicity questions, is that global sections
of these bundles vanishing on ample divisors provide algebraic differential
equations for any entire curve $f:\mathbb{C}\rightarrow X\backslash D$ where
$D$ is a normal crossing divisor on $X.$ Therefore, the main point is to
produce enough algebraically independent global holomorphic logarithmic jet
differentials. In the case of $\mathbb{P}^{3}\backslash X$ for a smooth
hypersurface $X\subset \mathbb{P}^{3},$ we have proved the existence of
global logarithmic jet differentials when $\deg (X)\geq 92$ in \cite{Rou05}.
Therefore to produce enough logarithmic jet differentials we take the
derivative of the logarithmic jet differential in the direction of the
vector fields constructed in the first part, just as in the compact case
\cite{Rou06}.

\section{Logarithmic jet bundles}

In this section we recall the basic facts about logarithmic jet bundles
following G. Dethloff and S. Lu \cite{DL96}.

Let $X$ be a complex manifold of dimension $n$. Let $x\in X.$ We consider
germs $f:(\mathbb{C},0)\rightarrow (X,x)$ of holomorphic curves. Then the
usual $k$-jet bundle, $J_{k}X,$ is the holomorphic fibre bundle whose fiber $%
J_{k}X_{x}$ is the set of equivalence classes of germs, $j_{k}(f),$ where
two germs are equivalent if they have the same Taylor expansions of order $%
k. $ Let $\pi :J_{k}X\rightarrow X$ be the natural projection.

Let $T_{X}^{\ast }$ be the holomorphic cotangent bundle over $X.$ Take a
holomorphic section $\omega \in H^{0}(O,T_{X}^{\ast })$ for some open subset
$O.$ For $j_{k}(f)\in J_{k}X_{\left| O\right. },$ we have $f^{\ast }\omega
=Z(t)dt$ and a well defined holomorphic mapping
\begin{equation*}
\widetilde{\omega }:J_{k}X_{\left| O\right. }\rightarrow \mathbb{C}%
^{k};j_{k}(f)\rightarrow \left( \frac{d^{j}Z}{dt^{j}}(0)\right) _{0\leq
j\leq k-1}.
\end{equation*}

If, moreover $\omega _{1},...,\omega _{n}$ are holomorphic 1-forms on $O$
such that $\omega _{1}\wedge ...\wedge \omega _{n}$ does not vanish
anywhere, then we have a biholomorphic map
\begin{equation*}
\left( \widetilde{\omega _{1}},...,\widetilde{\omega _{n}}\right) \times \pi
:J_{k}X_{\left| O\right. }\rightarrow \left( \mathbb{C}^{k}\right)
^{n}\times O
\end{equation*}

which gives the trivialization associated to $\omega _{1},...,\omega _{n}.$

Let $\overline{X}$ be a complex manifold with a normal crossing divisor $D.$
The pair $(\overline{X},D)$ is called a log manifold. Let $X=\overline{X}%
\backslash D.$

The logarithmic cotangent sheaf $\overline{T}_{X}^{\ast }=T_{\overline{X}%
}^{\ast }(\log D)$ is defined as the locally free subsheaf of the sheaf of
meromorphic 1-forms on $\overline{X},$ whose restriction to $X$ is $%
T_{X}^{\ast }$ and whose localization at any point $x\in D$ is given by
\begin{equation*}
\overline{T}_{X,x}^{\ast }=\underset{i=1}{\overset{l}{\sum }}\mathcal{O}_{%
\overline{X},x}\frac{dz_{i}}{z_{i}}+\underset{j=1+1}{\overset{n}{\sum }}%
\mathcal{O}_{\overline{X},x}dz_{j}
\end{equation*}
where the local coordinates $z_{1,}...,z_{n}$ around $x$ are chosen such
that $D=\{$ $z_{1}...z_{l}=0\}.$

Its dual, the logarithmic tangent sheaf $\overline{T}_{X}=T_{\overline{X}%
}(-\log D)$ is a locally free subsheaf of the holomorphic tangent bundle $T_{%
\overline{X}},$ whose restriction to $X$ is $T_{X}$ and whose localization
at any point $x\in D$ is given by
\begin{equation*}
\overline{T}_{X,x}=\underset{i=1}{\overset{l}{\sum }}\mathcal{O}_{\overline{X%
},x}z_{i}\frac{\partial }{\partial z_{i}}+\underset{j=1+1}{\overset{n}{\sum }%
}\mathcal{O}_{\overline{X},x}\frac{\partial }{\partial z_{j}}.
\end{equation*}

Given log-manifolds $(\overline{X},D)$ and $(\overline{X}^{\prime
},D^{\prime })$, a holomorphic map $F:\overline{X}^{\prime }\rightarrow
\overline{X}$ such that $F^{-1}(D)\subset D^{\prime }$ is called a
log-morphism from $(\overline{X}^{\prime },D^{\prime })$ to $(\overline{X}%
,D).$ It induces vector bundle morphisms
\begin{eqnarray*}
F^{\ast } &:&\overline{T}_{X}^{\ast }\rightarrow \overline{T}_{X^{\prime
}}^{\ast }; \\
F_{\ast } &:&\overline{T}_{X^{\prime }}\rightarrow \overline{T}_{X}.
\end{eqnarray*}

Let $s\in H^{0}(O,J_{k}\overline{X})$ be a holomorphic section over an open
subset $O\subset \overline{X}.$ We say that $s$ is a logarithmic $k$-jet
field if the map $\widetilde{\omega }\circ s_{\left| O^{\prime }\right.
}:O^{\prime }\rightarrow \mathbb{C}^{k}$ is holomorphic for all $\omega \in
H^{0}(O^{\prime },\overline{T}_{X}^{\ast })$ for all open subsets $O^{\prime
}$ of $O.$ The set of logarithmic $k$-jet fields over open subsets of $%
\overline{X}$ defines a subsheaf of the sheaf $J_{k}\overline{X},$ which we
denote by $\overline{J}_{k}X.$ $\overline{J}_{k}X$ is the sheaf of sections
of a holomorphic fibre bundle over $\overline{X},$ denoted again $\overline{J%
}_{k}X$ and called the logarithmic $k-$jet bundle of $(\overline{X},D).$

A log-morphism $F:(\overline{X}^{\prime },D^{\prime })\rightarrow (\overline{%
X},D)$ induces a canonical map
\begin{equation*}
F_{k}:\overline{J}_{k}X^{\prime }\rightarrow \overline{J}_{k}X.
\end{equation*}

We can express the local triviality of $\overline{J}_{k}X$ explicitly in
terms of coordinates. Let $z_{1,}...,z_{n}$ be coordinates in an open set $%
U\subset \overline{X}$ in which $D=\{z_{1}z_{2}...z_{l}=0\}.$ Let $\omega
_{1}=\frac{dz_{1}}{z_{1}},...,\omega _{l}=\frac{dz_{l}}{z_{l}},\omega
_{l+1}=dz_{l+1},...,\omega _{n}=dz_{n}.$ Then we have a biholomorphic map
\begin{equation*}
\left( \widetilde{\omega _{1}},...,\widetilde{\omega _{n}}\right) \times \pi
:\overline{J}_{k}X_{\left| U\right. }\rightarrow \left( \mathbb{C}%
^{k}\right) ^{n}\times U.
\end{equation*}

Let $s\in H^{0}(U,\overline{J}_{k}X)$ be given by $s(x)=(\xi
_{j}^{(i)}(x),x) $ in this trivialization where the indices $i$ correspond
to the orders of derivative. Then the same $s$ considered as an element of $%
H^{0}(U,J_{k}\overline{X})$ and trivialized by $\omega
_{1}=dz_{1},...,\omega _{n}=dz_{n}$ is given by $s(x)=(\widehat{\xi }%
_{j}^{(i)}(x),x)$ where
\begin{equation*}
\widehat{\xi }_{j}^{(i)}=\left\{
\begin{array}{c}
z_{i}(\xi _{j}^{(i)}+g_{i}(\xi _{j}^{(1)},...,\xi _{j}^{(i-1)})):j\leq l \\
\xi _{j}^{(i)}:j\geq l+1
\end{array}
\right.
\end{equation*}

The $g_{i}$ are polynomials in the variables $\xi _{j}^{(1)},...,\xi
_{j}^{(i-1)},$ obtained by expressing first the different components $\xi
_{j}^{(i)}$ of $\left( \widetilde{\frac{dz_{i}}{z_{i}}}\right) \circ s(x)$
in terms of the components $\widehat{\xi }_{j}^{(i)}$ of \ the components $%
\widehat{\xi }_{j}^{(i)}$ of $\widetilde{dz_{i}}\circ s(x)$ by using the
chain rule, and then by inverting this system.

\section{Logarithmic vector fields}

Let $\mathcal{X}\subset \mathbb{P}^{3}\times \mathbb{P}^{N_{d}}$ be the
universal surface of degree $d$ given by the equation
\begin{equation*}
\underset{\left| \alpha \right| =d}{\sum }a_{\alpha }Z^{\alpha }=0,\text{
where }[a]\in \mathbb{P}^{N_{d}}\text{ and }[Z]\in \mathbb{P}^{3}.
\end{equation*}

In this section we generalize the approach used in \cite{PaRou} (see
Proposition 11 of that article) and \cite{Rou06} to logarithmic jet bundles.
We use the notations: for $\alpha =(\alpha _{0},...,\alpha _{3})\in \mathbb{N%
}^{4},$ $\left| \alpha \right| =\sum_{i}\alpha _{i}$ and if $%
Z=(Z_{0},Z_{1},Z_{2},Z_{3})$ are homogeneous coordinates on $\mathbb{P}^{3},$
then $Z^{\alpha }=\prod Z_{j}^{\alpha _{j}}.$ $\mathcal{X}$ is a smooth
hypersurface of degree $(d,1)$ in $\mathbb{P}^{3}\times \mathbb{P}^{N_{d}}.$

We consider the log-manifold $(\mathbb{P}^{3}\times \mathbb{P}^{N_{d}},%
\mathcal{X)}$. We denote by $\overline{J_{3}}(\mathbb{P}^{3}\times \mathbb{P}%
^{N_{d}})$ the manifold of the logarithmic 3-jets$,$ and $\overline{J_{3}^{v}%
}(\mathbb{P}^{3}\times \mathbb{P}^{N_{d}})$ the submanifold of $\overline{%
J_{3}}(\mathbb{P}^{3}\times \mathbb{P}^{N_{d}})$ consisting of 3-jets
tangent to the fibers of the projection $\pi _{2}:\mathbb{P}^{3}\times
\mathbb{P}^{N_{d}}\rightarrow \mathbb{P}^{N_{d}}.$

We are going to construct meromorphic vector fields on $\overline{J_{3}^{v}}(%
\mathbb{P}^{3}\times \mathbb{P}^{N_{d}}).$

Let us consider
\begin{equation*}
\mathcal{Y}=(a_{d}Z_{4}^{d}+\underset{\left| \alpha \right| =d}{\sum }%
a_{\alpha }Z^{\alpha }=0)\subset \mathbb{P}^{4}\times U
\end{equation*}
where $U:=(a_{0...0d}\neq 0)\cap \left( \underset{\left| \alpha \right|
=d,\alpha _{n+2}=0}{\cup }(a_{\alpha }\neq 0)\right) \subset \mathbb{P}%
^{N_{d}+1}.$ We have the projection $\pi :\mathcal{Y}\rightarrow \mathbb{P}%
^{3}\times \mathbb{P}^{N_{d}}$ and $\pi ^{-1}(\mathcal{X})=(Z_{4}=0)$ $:=H$
therefore we obtain a log-morphism $\pi :(\mathcal{Y},H)\rightarrow (\mathbb{%
P}^{3}\times \mathbb{P}^{N_{d}},\mathcal{X})$ which induces a dominant map
\begin{equation*}
\pi _{3}:\overline{J_{3}^{v}}(\mathcal{Y})\rightarrow \overline{J_{3}^{v}}(%
\mathbb{P}^{3}\times \mathbb{P}^{N_{d}}).
\end{equation*}

Let us consider the set $\Omega _{0}:=(Z_{0}\neq 0)\times (a_{d}\neq
0)\subset \mathbb{P}^{4}\times U.$ We assume that global coordinates are
given on $\mathbb{C}^{4}$ and $\mathbb{C}^{N_{d}+1}.$ The equation of $%
\mathcal{Y}$ becomes
\begin{equation*}
\mathcal{Y}_{0}:=(z_{4}^{d}+\underset{\alpha }{\sum }a_{\alpha }z^{\alpha
}=0).
\end{equation*}

Following \cite{DL96} as explained above, we can obtain explicitly a
trivialization of $\overline{J_{3}}(\Omega _{0}).$ Let $\omega
^{1}=dz_{1},\omega ^{2}=dz_{2},\omega ^{3}=dz_{3},\omega ^{4}=\frac{dz_{4}}{%
z_{4}}.$ Then we have a biholomorphic map
\begin{equation*}
\overline{J_{3}}(\Omega _{0})\rightarrow \mathbb{C}^{4}\times U\times
\mathbb{C}^{4}\times \mathbb{C}^{4}\times \mathbb{C}^{4}
\end{equation*}
where the coordinates will be noted $(z_{i},a_{\alpha },\xi _{j}^{(i)}).$

Let's write the equations of $\overline{J_{3}^{v}}(\mathcal{Y}_{0})$ in this
trivialization. We have $\overline{J_{3}^{v}}(\mathcal{Y}_{0})=J_{3}^{v}(%
\mathcal{Y}_{0})\cap \overline{J_{3}}(\Omega _{0}).$ The equations of $%
J_{3}^{v}(\mathcal{Y}_{0})$ in the trivialization of $J_{3}(\Omega _{0})$
given by $\widehat{\omega }^{1}=dz_{1},\widehat{\omega }^{2}=dz_{2},\widehat{%
\omega }^{3}=dz_{3},\widehat{\omega }^{4}=dz_{4}\mathbb{\ }$can be written
in $\mathbb{C}^{4}\times U\times \mathbb{C}^{4}\times \mathbb{C}^{4}\times
\mathbb{C}^{4}$ with coordinates $(z_{i},a_{\alpha },\widehat{\xi }%
_{j}^{(i)})$:
\begin{equation*}
z_{4}^{d}+\underset{\left| \alpha \right| \leq d}{\sum }a_{\alpha }z^{\alpha
}=0
\end{equation*}
\begin{equation*}
dz_{4}^{d-1}\widehat{\xi }_{4}^{(1)}+\underset{j=1}{\overset{3}{\sum }}%
\underset{\left| \alpha \right| \leq d}{\sum }a_{\alpha }\frac{\partial
z^{\alpha }}{\partial z_{j}}\widehat{\xi }_{j}^{(1)}=0
\end{equation*}
\begin{equation*}
dz_{4}^{d-1}\widehat{\xi }_{4}^{(2)}+d(d-1)z_{4}^{d-2}\left( \widehat{\xi }%
_{4}^{(1)}\right) ^{2}+\underset{j=1}{\overset{3}{\sum }}\underset{\left|
\alpha \right| \leq d}{\sum }a_{\alpha }\frac{\partial z^{\alpha }}{\partial
z_{j}}\widehat{\xi }_{j}^{(2)}+\underset{j,k=1}{\overset{3}{\sum }}\underset{%
\left| \alpha \right| \leq d}{\sum }a_{\alpha }\frac{\partial ^{2}z^{\alpha }%
}{\partial z_{j}\partial z_{k}}\widehat{\xi }_{j}^{(1)}\widehat{\xi }%
_{k}^{(1)}=0
\end{equation*}
\begin{eqnarray*}
dz_{4}^{d-1}\widehat{\xi }_{4}^{(3)}+3d(d-1)z_{4}^{d-2}\widehat{\xi }%
_{4}^{(1)}\widehat{\xi }_{4}^{(2)}+d(d-1)(d-2)z_{4}^{d-3}\left( \widehat{\xi
}_{4}^{(1)}\right) ^{3} +\underset{j=1}{\overset{3}{\sum }}\underset{\left|
\alpha \right| \leq d}{\sum }a_{\alpha }\frac{\partial z^{\alpha }}{\partial
z_{j}}\widehat{\xi }_{j}^{(3)}&&  \notag \\
+\underset{j,k=1}{\overset{3}{3\sum }}\underset{\left| \alpha \right| \leq d%
}{\sum }a_{\alpha }\frac{\partial ^{2}z^{\alpha }}{\partial z_{j}\partial
z_{k}}\widehat{\xi }_{j}^{(2)}\widehat{\xi }_{k}^{(1)} +\underset{j,k,l=1}{%
\overset{3}{\sum }}\underset{\left| \alpha \right| \leq d}{\sum }a_{\alpha }%
\frac{\partial ^{3}z^{\alpha }}{\partial z_{j}\partial z_{k}\partial z_{l}}%
\widehat{\xi }_{j}^{(1)}\widehat{\xi }_{k}^{(1)}\widehat{\xi }_{l}^{(1)}=0&&
\end{eqnarray*}

The relations between the two systems of coordinates can be computed as
explained above and are given by
\begin{eqnarray*}
\widehat{\xi }_{j}^{(i)} &=&\xi _{j}^{(i)}\text{ for }j\leq 3 \\
\widehat{\xi }_{4}^{(1)} &=&z_{4}\xi _{4}^{(1)} \\
\widehat{\xi }_{4}^{(2)} &=&z_{4}(\xi _{4}^{(2)}+\left( \xi
_{4}^{(1)}\right) ^{2}) \\
\widehat{\xi }_{4}^{(3)} &=&z_{4}(\xi _{4}^{(3)}+3\xi _{4}^{(1)}\xi
_{4}^{(2)}+\left( \xi _{4}^{(1)}\right) ^{3})
\end{eqnarray*}

Therefore, to obtain the equations of $\overline{J_{3}^{v}}(\mathcal{Y}_{0})$
in the first trivialization, we just have to substitute the previous
relations
\begin{equation}
z_{4}^{d}+\underset{\left| \alpha \right| \leq d}{\sum }a_{\alpha }z^{\alpha
}=0
\end{equation}

\begin{equation}
dz_{4}^{d}\xi _{4}^{(1)}+\underset{j=1}{\overset{3}{\sum }}\underset{\left|
\alpha \right| \leq d}{\sum }a_{\alpha }\frac{\partial z^{\alpha }}{\partial
z_{j}}\xi _{j}^{(1)}=0
\end{equation}

\begin{equation}
dz_{4}^{d}\xi _{4}^{(2)}+d^{2}z_{4}^{d}\left( \xi _{4}^{(1)}\right) ^{2}+%
\underset{j=1}{\overset{3}{\sum }}\underset{\left| \alpha \right| \leq d}{%
\sum }a_{\alpha }\frac{\partial z^{\alpha }}{\partial z_{j}}\xi _{j}^{(2)}+%
\underset{j,k=1}{\overset{3}{\sum }}\underset{\left| \alpha \right| \leq d}{%
\sum }a_{\alpha }\frac{\partial ^{2}z^{\alpha }}{\partial z_{j}\partial z_{k}%
}\xi _{j}^{(1)}\xi _{k}^{(1)}=0
\end{equation}

\begin{eqnarray}
dz_{4}^{d}\xi _{4}^{(3)}+3d^{2}z_{4}^{d}\xi _{4}^{(1)}\xi
_{4}^{(2)}+d^{3}z_{4}^{d}\left( \xi _{4}^{(1)}\right) ^{3}+\underset{j=1}{%
\overset{3}{\sum }}\underset{\left| \alpha \right| \leq d}{\sum }a_{\alpha }%
\frac{\partial z^{\alpha }}{\partial z_{j}}\xi _{j}^{(3)} &&  \notag \\
+3\underset{j,k=1}{\overset{3}{\sum }}\underset{\left| \alpha \right| \leq d%
}{\sum }a_{\alpha }\frac{\partial ^{2}z^{\alpha }}{\partial z_{j}\partial
z_{k}}\xi _{j}^{(2)}\xi _{k}^{(1)}+\underset{j,k,l=1}{\overset{3}{\sum }}%
\underset{\left| \alpha \right| \leq d}{\sum }a_{\alpha }\frac{\partial
^{3}z^{\alpha }}{\partial z_{j}\partial z_{k}\partial z_{l}}\xi
_{j}^{(1)}\xi _{k}^{(1)}\xi _{l}^{(1)}=0 &&
\end{eqnarray}

Following the method used in \cite{Rou06} for the compact case, we are going
to prove that $T_{\overline{J_{3}^{v}}(\mathcal{Y})}\otimes \mathcal{O}_{%
\mathbb{P}^{4}}(c)\otimes \mathcal{O}_{\mathbb{P}^{N_{d}}+1}(\ast )$ is
generated by its global sections on $\overline{J_{3}^{v}}(\mathcal{Y}%
)\backslash (\Sigma \cup p^{-1}(H)),$ where $p:$ $\overline{J_{3}^{v}}(%
\mathcal{Y})\rightarrow $ $\mathcal{Y}$ is the natural projection, $\Sigma $
a subvariety that will be defined below, and $c\in \mathbb{N}$ a constant
independant of $d.$ Consider a vector field
\begin{equation*}
V=\underset{\left| \alpha \right| \leq d}{\sum }v_{\alpha }\frac{\partial }{%
\partial a_{\alpha }}+\underset{j}{\sum }v_{j}\frac{\partial }{\partial z_{j}%
}+\underset{j,k}{\sum }w_{j}^{(k)}\frac{\partial }{\partial \xi _{j}^{(k)}}
\end{equation*}
on $\mathbb{C}^{4}\times U\times \mathbb{C}^{4}\times
\mathbb{C}^{4}\times
\mathbb{C}^{4}.$ The conditions to be satisfied by $V$ to be tangent to $%
\overline{J_{3}^{v}}(\mathcal{Y}_{0})$ are the following
\begin{equation}
\underset{\left| \alpha \right| \leq d}{\sum }v_{\alpha }z^{\alpha }+%
\underset{j=1}{\overset{3}{\sum }}\underset{\left| \alpha \right| \leq d}{%
\sum }a_{\alpha }\frac{\partial z^{\alpha }}{\partial z_{j}}%
v_{j}+dz_{4}^{d-1}v_{4}=0
\end{equation}

\begin{eqnarray}
\underset{j=1}{\overset{3}{\sum }}\underset{\left| \alpha \right| \leq d,%
\text{ }\alpha _{1}<d}{\sum }v_{\alpha }\frac{\partial z^{\alpha }}{\partial
z_{j}}\xi _{j}^{(1)}+\underset{j,k=1}{\overset{3}{\sum }}\underset{\left|
\alpha \right| \leq d}{\sum }a_{\alpha }\frac{\partial ^{2}z^{\alpha }}{%
\partial z_{j}\partial z_{k}}v_{j}\xi _{k}^{(1)}+\underset{j=1}{\overset{3}{%
\sum }}\underset{\left| \alpha \right| \leq d}{\sum }a_{\alpha }\frac{%
\partial z^{\alpha }}{\partial z_{j}}w_{j}^{(1)} &&  \notag \\
+d^{2}z_{4}^{d-1}v_{4}\xi _{4}^{(1)}+dz_{4}^{d}w_{4}^{(1)}=0 &&
\end{eqnarray}
\begin{eqnarray}
\underset{\left| \alpha \right| \leq d}{\sum }(\underset{j=1}{\overset{3}{%
\sum }}\frac{\partial z^{\alpha }}{\partial z_{j}}\xi _{j}^{(2)}+\underset{%
j,k=1}{\overset{3}{\sum }}\frac{\partial ^{2}z^{\alpha }}{\partial
z_{j}\partial z_{k}}\xi _{j}^{(1)}\xi _{k}^{(1)})v_{\alpha } &&  \notag \\
+\underset{j=1}{\overset{3}{\sum }}\underset{\left| \alpha \right| \leq d}{%
\sum }a_{\alpha }(\underset{k=1}{\overset{3}{\sum }}\frac{\partial
^{2}z^{\alpha }}{\partial z_{j}\partial z_{k}}\xi _{k}^{(2)}+\underset{k,l=1%
}{\overset{3}{\sum }}\frac{\partial ^{3}z^{\alpha }}{\partial z_{j}\partial
z_{k}\partial z_{l}}\xi _{k}^{(1)}\xi _{l}^{(1)})v_{j} &&  \notag \\
+\underset{\left| \alpha \right| \leq d}{\sum }(\underset{j,k=1}{\overset{3}{%
\sum }}a_{\alpha }\frac{\partial ^{2}z^{\alpha }}{\partial z_{j}\partial
z_{k}}(w_{j}^{(1)}\xi _{k}^{(1)}+w_{k}^{(1)}\xi _{j}^{(1)})+\underset{j=1}{%
\overset{3}{\sum }}a_{\alpha }\frac{\partial z^{\alpha }}{\partial z_{j}}%
w_{j}^{(2)}) &&  \notag \\
+v_{4}d^{2}z_{4}^{d-1}(\xi _{4}^{(2)}+d\left( \xi _{4}^{(1)}\right)
^{2})+2d^{2}z_{4}^{d}w_{4}^{(1)}\xi _{4}^{(1)}+dz_{4}^{d}w_{4}^{(2)}=0 &&
\end{eqnarray}
\begin{eqnarray}
\underset{\left| \alpha \right| \leq d}{\sum }(\underset{j=1}{\overset{3}{%
\sum }}\frac{\partial z^{\alpha }}{\partial z_{j}}\xi _{j}^{(3)}+3\underset{%
j,k=1}{\overset{3}{\sum }}\frac{\partial ^{2}z^{\alpha }}{\partial
z_{j}\partial z_{k}}\xi _{j}^{(2)}\xi _{k}^{(1)}+\underset{j,k,l=1}{\overset{%
3}{\sum }}\frac{\partial ^{3}z^{\alpha }}{\partial z_{j}\partial
z_{k}\partial z_{l}}\xi _{j}^{(1)}\xi _{k}^{(1)}\xi _{l}^{(1)})v_{\alpha } &&
\notag \\
+\underset{j=1}{\overset{3}{\sum }}\underset{\left| \alpha \right| \leq d}{%
\sum }a_{\alpha }(\underset{k=1}{\overset{3}{\sum }}\frac{\partial
^{2}z^{\alpha }}{\partial z_{j}\partial z_{k}}\xi _{k}^{(3)}+3\underset{k,l=1%
}{\overset{3}{\sum }}\frac{\partial ^{3}z^{\alpha }}{\partial z_{j}\partial
z_{k}\partial z_{l}}\xi _{k}^{(2)}\xi _{l}^{(1)}&&  \notag \\
+\underset{k,l,m=1}{\overset{%
3}{\sum }}\frac{\partial ^{4}z^{\alpha }}{\partial z_{j}\partial
z_{k}\partial z_{l}\partial z_{m}}\xi _{k}^{(1)}\xi _{l}^{(1)}\xi
_{m}^{(1)})v_{j} &&  \notag \\
+\underset{\left| \alpha \right| \leq d}{\sum }(\underset{j,k,l=1}{\overset{3%
}{\sum }}a_{\alpha }\frac{\partial ^{3}z^{\alpha }}{\partial z_{j}\partial
z_{k}\partial z_{l}}(w_{j}^{(1)}\xi _{k}^{(1)}\xi _{l}^{(1)}+\xi
_{j}^{(1)}w_{k}^{(1)}\xi _{l}^{(1)}+\xi _{j}^{(1)}\xi _{k}^{(1)}w_{l}^{(1)})
&&  \notag \\
+3\underset{j,k=1}{\overset{3}{\sum }}a_{\alpha }\frac{\partial
^{2}z^{\alpha }}{\partial z_{j}\partial z_{k}}(w_{j}^{(2)}\xi _{k}^{(1)}+\xi
_{j}^{(2)}w_{k}^{(1)})+\underset{j=1}{\overset{3}{\sum }}a_{\alpha }\frac{%
\partial z^{\alpha }}{\partial z_{j}}w_{j}^{(3)})&&  \notag \\
+d^{2}z_{4}^{d-1}v_{4}(\xi
_{4}^{(3)}+3d\xi _{4}^{(1)}\xi _{4}^{(2)}+d\left( \xi
_{4}^{(1)}\right)
^{3}) &&  \notag \\
+dz_{4}^{d}w_{4}^{(3)}+3d^{2}z_{4}^{d}(\xi _{4}^{(2)}w_{4}^{(1)}+\xi
_{4}^{(1)}w_{4}^{(2)})+3d^{3}z_{4}^{d}w_{4}^{(1)}\left( \xi
_{4}^{(1)}\right) ^{2}=0 &&
\end{eqnarray}

We can introduce the first package of vector fields tangent to $\overline{%
J_{3}^{v}}(\mathcal{Y}_{0}).$ We denote by $\delta _{j}\in \mathbb{N}^{3}$
the multi-index whose j-component is equal to 1 and the other are zero.

For $\alpha _{1}\geq 4:$%
\begin{equation*}
V_{\alpha }^{400}:=\frac{\partial }{\partial a_{\alpha }}-4z_{1}\frac{%
\partial }{\partial a_{\alpha -\delta _{1}}}+6z_{1}^{2}\frac{\partial }{%
\partial a_{\alpha -2\delta _{1}}}-4z_{1}^{3}\frac{\partial }{\partial
a_{\alpha -3\delta _{1}}}+z_{1}^{4}\frac{\partial }{\partial a_{\alpha
-4\delta _{1}}}.
\end{equation*}

For $\alpha _{1}\geq 3,\alpha _{2}\geq 1:$%
\begin{eqnarray*}
V_{\alpha }^{310} &:&=\frac{\partial }{\partial a_{\alpha }}-3z_{1}\frac{%
\partial }{\partial a_{\alpha -\delta _{1}}}-z_{2}\frac{\partial }{\partial
a_{\alpha -\delta _{2}}}+3z_{1}z_{2}\frac{\partial }{\partial a_{\alpha
-\delta _{1}-\delta _{2}}} \\
&&+3z_{1}^{2}\frac{\partial }{\partial a_{\alpha -2\delta _{1}}}%
-3z_{1}^{2}z_{2}\frac{\partial }{\partial a_{\alpha -2\delta _{1}-\delta
_{2}}}-z_{1}^{3}\frac{\partial }{\partial a_{\alpha -3\delta _{1}}}%
+z_{1}^{3}z_{2}\frac{\partial }{\partial a_{\alpha -3\delta _{1}-\delta _{2}}%
}.
\end{eqnarray*}

For $\alpha _{1}\geq 2,\alpha _{2}\geq 2:$%
\begin{eqnarray*}
V_{\alpha }^{220}&:&=\frac{\partial }{\partial a_{\alpha }}-z_{2}\frac{%
\partial }{\partial a_{\alpha -\delta _{2}}}-z_{1}\frac{\partial }{\partial
a_{\alpha -\delta _{1}}}+z_{1}z_{2}^{2}\frac{\partial }{\partial
a_{\alpha -\delta _{1}-2\delta
_{2}}}\\
&&+z_{1}^{2}z_{2}\frac{\partial }{\partial
a_{\alpha -2\delta _{1}-\delta _{2}}}-z_{1}^{2}z_{2}^{2}\frac{\partial }{%
\partial a_{\alpha -2\delta _{1}-2\delta _{2}}}.
\end{eqnarray*}

For $\alpha _{1}\geq 2,\alpha _{2}\geq 1,\alpha _{3}\geq 1:$%
\begin{eqnarray*}
V_{\alpha }^{211} &:&=\frac{\partial }{\partial a_{\alpha }}-z_{3}\frac{%
\partial }{\partial a_{\alpha -\delta _{3}}}-z_{2}\frac{\partial }{\partial
a_{\alpha -\delta _{2}}}-2z_{1}\frac{\partial }{\partial a_{\alpha -\delta
_{1}}}+z_{2}z_{3}\frac{\partial }{\partial a_{\alpha -\delta _{2}-\delta
_{3}}} \\
&&+2z_{1}z_{3}\frac{\partial }{\partial a_{\alpha -\delta _{1}-\delta _{3}}}%
+2z_{1}z_{2}\frac{\partial }{\partial a_{\alpha -\delta _{1}-\delta _{2}}}%
+z_{1}^{2}\frac{\partial }{\partial a_{\alpha -2\delta _{1}}} \\
&&-2z_{1}z_{2}z_{3}\frac{\partial }{\partial a_{\alpha -\delta
_{1}-\delta _{2}-\delta _{3}}}-z_{1}^{2}z_{3}\frac{\partial
}{\partial a_{\alpha -2\delta _{1}-\delta
_{3}}}\\
&&-z_{1}^{2}z_{2}\frac{\partial }{\partial
a_{\alpha -2\delta _{1}-\delta _{2}}}+z_{1}^{2}z_{2}z_{3}\frac{\partial }{%
\partial a_{\alpha -2\delta _{1}-\delta _{2}-\delta _{3}}}.
\end{eqnarray*}

Similar vector fields are constructed by permuting the z-variables, and
changing the index $\alpha $ as indicated by the permutation. The pole order
of the previous vector fields is equal to 4.

\begin{lemma}
For any $(v_{i})_{1\leq i\leq 4}\in \mathbb{C}^{4},$ there exist $v_{\alpha
}(a),$ with degree at most 1 in the variables $(a_{\gamma }),$ such that $V:=%
\underset{\alpha }{\sum }v_{\alpha }(a)\frac{\partial }{\partial a_{\alpha }}%
+\underset{1\leq j\leq 3}{\sum }v_{j}\frac{\partial }{\partial z_{j}}%
+v_{4}z_{4}\frac{\partial }{\partial z_{4}}$ is tangent to $\overline{%
J_{3}^{v}}(\mathcal{Y}_{0})$ at each point.
\end{lemma}

\begin{proof}
First, we substitute equations 1, 2, 3, 4 in equations 5, 6, 7, 8 to get rid
of $z_{4},\xi _{4}^{(i)}(1\leq i\leq 3).$ Then, we impose the additional
conditions of vanishing for the coefficients of $\xi _{j}^{(1)}$ in the
second equation (respectively of $\xi _{j}^{(1)}\xi _{k}^{(1)}$ in the third
equation and $\xi _{j}^{(1)}\xi _{k}^{(1)}\xi _{l}^{(1)}$ in the fourth
equation) for any $1\leq j\leq k\leq l\leq 3$. Then the coefficients of $\xi
_{j}^{(2)}$ (respectively $\xi _{j}^{(2)}\xi _{k}^{(1)}$ and $\xi
_{j}^{(3)}) $ are automatically zero in the third (respectively fourth)
equation. The resulting equations are
\begin{equation*}
\underset{\left| \alpha \right| \leq d}{\sum }v_{\alpha }z^{\alpha }+%
\underset{j=1}{\overset{3}{\sum }}\underset{\left| \alpha \right| \leq d}{%
\sum }a_{\alpha }\frac{\partial z^{\alpha }}{\partial z_{j}}v_{j}-dv_{4}%
\underset{\left| \alpha \right| \leq d}{\sum }a_{\alpha }z^{\alpha }=0
\end{equation*}

\begin{equation*}
\underset{\left| \alpha \right| \leq d}{\sum }v_{\alpha }\frac{\partial
z^{\alpha }}{\partial z_{j}}+\underset{k=1}{\overset{3}{\sum }}\underset{%
\left| \alpha \right| \leq d}{\sum }a_{\alpha }\frac{\partial ^{2}z^{\alpha }%
}{\partial z_{j}\partial z_{k}}v_{k}-dv_{4}\underset{\left| \alpha \right|
\leq d}{\sum }a_{\alpha }\frac{\partial z^{\alpha }}{\partial z_{j}}=0
\end{equation*}

\begin{equation*}
\underset{\left| \alpha \right| \leq d}{\sum }\frac{\partial ^{2}z^{\alpha }%
}{\partial z_{j}\partial z_{k}}v_{\alpha }+\underset{l=1}{\overset{3}{\sum }}%
\underset{\left| \alpha \right| \leq d}{\sum }a_{\alpha }\frac{\partial
^{3}z^{\alpha }}{\partial z_{j}\partial z_{k}\partial z_{l}}v_{l}-dv_{4}%
\underset{\left| \alpha \right| \leq d}{\sum }a_{\alpha }\frac{\partial
^{2}z^{\alpha }}{\partial z_{j}\partial z_{k}}=0
\end{equation*}

\begin{equation*}
\underset{\left| \alpha \right| \leq d}{\sum }\frac{\partial ^{3}z^{\alpha }%
}{\partial z_{j}\partial z_{k}\partial z_{l}}v_{\alpha }+\underset{m=1}{%
\overset{3}{\sum }}\underset{\left| \alpha \right| \leq d}{\sum }a_{\alpha }%
\frac{\partial ^{4}z^{\alpha }}{\partial z_{j}\partial z_{k}\partial
z_{l}\partial z_{m}}v_{m}-dv_{4}\underset{\left| \alpha \right| \leq d}{\sum
}a_{\alpha }\frac{\partial ^{3}z^{\alpha }}{\partial z_{j}\partial
z_{k}\partial z_{l}}=0
\end{equation*}

Now we can observe that if the $v_{\alpha }(a)$ satisfy the first equation,
they automatically satisfy the other ones because the $v_{\alpha }$ are
constants with respect to $z$. Therefore it is sufficient to find $%
(v_{\alpha })$ satisfying the first equation. We identify the coefficients
of $z^{\rho }=z_{1}^{\rho _{1}}$ $z_{2}^{\rho _{2}}$ $z_{3}^{\rho _{3}}:$%
\begin{equation*}
v_{\rho }+\underset{j=1}{\overset{4}{\sum }}a_{\rho +\delta _{j}}v_{j}(\rho
_{j}+1)-dv_{4}a_{\rho }=0.
\end{equation*}
\end{proof}

Another family of vector fields can be obtained in the following way.
Consider a $4\times 4$-matrix $A=\left(
\begin{array}{cccc}
A_{1}^{1} & A_{1}^{2} & A_{1}^{3} & 0 \\
A_{2}^{1} & A_{2}^{2} & A_{2}^{3} & 0 \\
A_{3}^{1} & A_{3}^{2} & A_{3}^{3} & 0 \\
A_{4}^{1} & A_{4}^{2} & A_{4}^{3} & 0
\end{array}
\right) \in \mathcal{M}_{4}(\mathbb{C})$ and let $\widetilde{V}:=\underset{%
j,k}{\sum }w_{j}^{(k)}\frac{\partial }{\partial \xi _{j}^{(k)}},$ where $%
w^{(k)}:=A\xi ^{(k)},$ for $k=1,2,3.$

\begin{lemma}
There exist polynomials $v_{\alpha }(z,a):=\underset{\left| \beta \right|
\leq 3}{\sum }v_{\beta }^{\alpha }(a)z^{\beta }$ where each coefficient $%
v_{\beta }^{\alpha }$ has degree at most 1 in the variables $(a_{\gamma })$
such that
\begin{equation*}
V:=\underset{\alpha }{\sum }v_{\alpha }(z,a)\frac{\partial }{\partial
a_{\alpha }}+\widetilde{V}
\end{equation*}
is tangent to $\overline{J_{3}^{v}}(\mathcal{Y}_{0})$ at each point.
\end{lemma}

\begin{proof}
First, we substitute equations 1, 2, 3, 4 in equations 5, 6, 7, 8 to get rid
of $z_{4},\xi _{4}^{(i)}(1\leq i\leq 3).$ We impose the additional
conditions of vanishing for the coefficients of $\xi _{j}^{(1)}$ in the
second equation (respectively of $\xi _{j}^{(1)}\xi _{k}^{(1)}$ in the third
equation and $\xi _{j}^{(1)}\xi _{k}^{(1)}\xi _{l}^{(1)}$ in the fourth
equation) for any $1\leq j\leq k\leq l\leq 3$. Then the coefficients of $\xi
_{j}^{(2)}$ (respectively $\xi _{j}^{(2)}\xi _{k}^{(1)}$ and $\xi
_{j}^{(3)}) $ are automatically zero in the third (respectively fourth)
equation. The resulting equations are
\begin{equation*}
\underset{\left| \alpha \right| \leq d}{\sum }v_{\alpha }z^{\alpha }=0\text{
\ }(9)
\end{equation*}

\begin{equation*}
\underset{\left| \alpha \right| \leq d}{\sum }v_{\alpha }\frac{\partial
z^{\alpha }}{\partial z_{j}}+\underset{k=1}{\overset{3}{\sum }}\underset{%
\left| \alpha \right| \leq d}{\sum }a_{\alpha }\frac{\partial z^{\alpha }}{%
\partial z_{k}}A_{k}^{j}-dA_{4}^{j}\underset{\left| \alpha \right| \leq d}{%
\sum }a_{\alpha }z^{\alpha }=0\text{ \ \ }(10_{j})
\end{equation*}

\begin{equation*}
\underset{\alpha }{\sum }\frac{\partial ^{2}z^{\alpha }}{\partial
z_{j}\partial z_{k}}v_{\alpha }+\underset{\alpha ,p}{\sum }a_{\alpha }\frac{%
\partial ^{2}z^{\alpha }}{\partial z_{j}\partial z_{p}}A_{p}^{k}+\underset{%
\alpha ,p}{\sum }a_{\alpha }\frac{\partial ^{2}z^{\alpha }}{\partial
z_{k}\partial z_{p}}A_{p}^{j}-2dA_{4}^{j}\underset{\left| \alpha \right|
\leq d}{\sum }a_{\alpha }\frac{\partial z^{\alpha }}{\partial z_{k}}=0\text{
\ }(11_{jk})
\end{equation*}

\begin{eqnarray*}
\underset{\alpha }{\sum }\frac{\partial ^{3}z^{\alpha }}{\partial
z_{j}\partial z_{k}\partial z_{l}}v_{\alpha }+\underset{\alpha ,p}{\sum }%
a_{\alpha }\frac{\partial ^{3}z^{\alpha }}{\partial z_{p}\partial
z_{k}\partial z_{l}}A_{p}^{j}+\underset{\alpha ,p}{\sum }a_{\alpha }\frac{%
\partial ^{3}z^{\alpha }}{\partial z_{j}\partial z_{p}\partial z_{l}}%
A_{p}^{k} && \notag \\
+\underset{\alpha ,p}{\sum }a_{\alpha
}\frac{\partial ^{3}z^{\alpha
}}{\partial z_{j}\partial z_{k}\partial z_{p}}A_{p}^{l}-3dA_{4}^{l}\underset{%
\left| \alpha \right| \leq d}{\sum }a_{\alpha }\frac{\partial ^{2}z^{\alpha }%
}{\partial z_{j}\partial z_{k}}=0\text{ \ }(12_{jkl})&&
\end{eqnarray*}

The equations for the unknowns $v_{\beta }^{\alpha }$ are obtained by
identifying the coefficients of the monomials $z^{\rho }$ in the above
equations.

The monomials $z^{\rho }$ in (9) are $z_{1}^{\rho _{1}}$
$z_{2}^{\rho _{2}}$ $z_{3}^{\rho _{3}}$ with $\sum \rho _{i}\leq
d$.

If all the components of $\rho $ are greater than 3, then we obtain the
following system

13. The coefficient of $z^{\rho }$ in (9) impose the condition
\begin{equation*}
\underset{\alpha +\beta =\rho }{\sum }v_{\beta }^{\alpha }=0
\end{equation*}

14$_{j}.$ The coefficient of the monomial $z^{\rho -\delta _{j}}$ in $%
(10_{j}) $ impose the condition
\begin{equation*}
\underset{\alpha +\beta =\rho }{\sum }\alpha _{j}v_{\beta }^{\alpha
}=l_{j}(a)
\end{equation*}
where $l_{j}$ is a linear expression in the $a$-variables.

14$_{jj}.$ For $j=1,...,3$ the coefficient of the monomial $z^{\rho -2\delta
_{j}}$ in $(11_{jj})$ impose the condition
\begin{equation*}
\underset{\alpha +\beta =\rho }{\sum }\alpha _{j}(\alpha _{j}-1)v_{\beta
}^{\alpha }=l_{jj}(a)
\end{equation*}

14$_{jk}.$ For $1\leq j<k\leq 3$ the coefficient of the monomial $z^{\rho
-\delta _{j}-\delta _{k}}$ in $(11_{jk})$ impose the condition
\begin{equation*}
\underset{\alpha +\beta =\rho }{\sum }\alpha _{j}\alpha _{k}v_{\beta
}^{\alpha }=l_{jk}(a)
\end{equation*}

15$_{jjj}.$ For $j=1,...,3$ the coefficient of the monomial $z^{\rho
-3\delta _{j}}$ in $(12_{jjj})$ impose the condition
\begin{equation*}
\underset{\alpha +\beta =\rho }{\sum }\alpha _{j}(\alpha _{j}-1)(\alpha
_{j}-2)v_{\beta }^{\alpha }=l_{jjj}(a)
\end{equation*}

15$_{jjk}.$ For $1\leq j<k\leq 3$ the coefficient of the monomial $z^{\rho
-2\delta _{j}-\delta _{k}}$ in $(12_{jjk})$ impose the condition
\begin{equation*}
\underset{\alpha +\beta =\rho }{\sum }\alpha _{j}(\alpha _{j}-1)\alpha
_{k}v_{\beta }^{\alpha }=l_{jjk}(a)
\end{equation*}

15$_{jkl}.$ For $1\leq j<k<l\leq 3$ the coefficient of the monomial $z^{\rho
-\delta _{j}-\delta _{k}-\delta _{l}}$ in $(12_{jjk})$ impose the condition
\begin{equation*}
\underset{\alpha +\beta =\rho }{\sum }\alpha _{j}\alpha _{k}\alpha
_{l}v_{\beta }^{\alpha }=l_{jkl}(a)
\end{equation*}

The determinant of the matrix associated to the system is not zero. Indeed,
for each $\rho $ the matrix whose column $C_{\beta }$ consists of the
partial derivatives of order at most 3 of the monomial $z^{\rho -\beta }$
has the same determinant, at the point $z_{0}=(1,1,1),$ as our system.
Therefore if the determinant is zero, we would have a non-identically zero
polynomial
\begin{equation*}
Q(z)=\underset{\beta }{\sum }a_{\beta }z^{\rho -\beta }
\end{equation*}
such that all its partial derivatives of order less or equal to 3 vanish at $%
z_{0}.$ Thus the same is true for
\begin{equation*}
P(z)=z^{\rho }Q(\frac{1}{z_{1}},...,\frac{1}{z_{3}})=\underset{\beta }{\sum }%
a_{\beta }z^{\beta }.
\end{equation*}
But this implies $P\equiv 0.$

Finally, we conclude by Cramer's rule. The systems we have to solve are
never over determined. The lemma is proved.
\end{proof}

\begin{remark}
We have chosen the matrix $A$ with this form because we are interested to
prove the global generation statement on $\overline{%
J_{3}^{v}}(\mathcal{Y})\backslash (\Sigma \cup p^{-1}(H))$ where $\Sigma $ is the closure of $%
\Sigma _{0}=\{(z,a,\xi ^{(1)},\xi ^{(2)},\xi ^{(3)})\in \overline{J_{3}^{v}}(%
\mathcal{Y}_{0})$ $/$ $\det \left( \xi _{i}^{(j)}\right) _{1\leq i,j\leq
3}=0\}$
\end{remark}

\begin{proposition}
\label{p1}The vector space $T_{\overline{J_{3}^{v}}(\mathcal{Y})}\otimes
\mathcal{O}_{\mathbb{P}^{4}}(12)\otimes \mathcal{O}_{\mathbb{P}%
^{N_{d}}+1}(\ast )$ is generated by its global sections on $\overline{%
J_{3}^{v}}(\mathcal{Y})\backslash (\Sigma \cup p^{-1}(H)).$
\end{proposition}

\begin{proof}
From the preceding lemmas, we are reduced to consider $V=\underset{\left|
\alpha \right| \leq 3}{\sum }v_{\alpha }\frac{\partial }{\partial a_{\alpha }%
}.$ The conditions for $V$ to be tangent to $\overline{J_{3}^{v}}(\mathcal{Y}%
_{0})$ are
\begin{equation*}
\underset{\left| \alpha \right| \leq 3}{\sum }v_{\alpha }z^{\alpha }=0
\end{equation*}

\begin{equation*}
\underset{j=1}{\overset{3}{\sum }}\underset{\left| \alpha \right| \leq 3}{%
\sum }v_{\alpha }\frac{\partial z^{\alpha }}{\partial z_{j}}\xi _{j}^{(1)}=0
\end{equation*}

\begin{equation*}
\underset{\left| \alpha \right| \leq 3}{\sum }(\underset{j=1}{\overset{3}{%
\sum }}\frac{\partial z^{\alpha }}{\partial z_{j}}\xi _{j}^{(2)}+\underset{%
j,k=1}{\overset{3}{\sum }}\frac{\partial ^{2}z^{\alpha }}{\partial
z_{j}\partial z_{k}}\xi _{j}^{(1)}\xi _{k}^{(1)})v_{\alpha }=0
\end{equation*}

\begin{equation*}
\underset{\left| \alpha \right| \leq 3}{\sum }(\underset{j=1}{\overset{3}{%
\sum }}\frac{\partial z^{\alpha }}{\partial z_{j}}\xi _{j}^{(3)}+3\underset{%
j,k=1}{\overset{3}{\sum }}\frac{\partial ^{2}z^{\alpha }}{\partial
z_{j}\partial z_{k}}\xi _{j}^{(2)}\xi _{k}^{(1)}+\underset{j,k,l=1}{\overset{%
3}{\sum }}\frac{\partial ^{3}z^{\alpha }}{\partial z_{j}\partial
z_{k}\partial z_{l}}\xi _{j}^{(1)}\xi _{k}^{(1)}\xi _{l}^{(1)})v_{\alpha }=0
\end{equation*}

We denote by $W_{jkl}$ the wronskian operator corresponding to the variables
$z_{j},z_{k},z_{l}.$ We have $W_{123}:=\det (\xi _{j}^{(i)})_{1\leq i,j\leq
3}\neq 0.$ Then we can solve the previous system with $%
v_{000},v_{100},v_{010},v_{001}$ as unknowns. By the Cramer rule, each of
the previous quantity is a linear combination of the $v_{\alpha },$ $\left|
\alpha \right| \leq 3,$ $\alpha \neq (000),$ $(100),$ $(010),$ $(001)$ with
coefficients rational functions in $z,\xi ^{(1)},\xi ^{(2)},\xi ^{(3)}.$ The
denominator is $W_{123}$ and the numerator is a polynomial whose monomials
verify either:

i) degree in $z$ at most 3 and degree in each $\xi ^{(i)}$ at most 1.

ii) degree in $z$ at most 2 and degree in $\xi ^{(1)}$ at most 3, degree in $%
\xi ^{(2)}$ at most 0, degree in $\xi ^{(3)}$ at most 1.

iii) degree in $z$ at most 2 and degree in $\xi ^{(1)}$ at most 2, degree in
$\xi ^{(2)}$ at most 2, degree in $\xi ^{(3)}$ at most 0.

iv) degree in $z$ at most 1 and degree in $\xi ^{(1)}$ at most 4, degree in $%
\xi ^{(2)}$ at most 1, degree in $\xi ^{(3)}$ at most 0.

$\xi ^{(1)}$ has a pole of order 2, $\xi ^{(2)}$ has a pole of order 3 and $%
\xi ^{(3)}$ has a pole of order 4, therefore the previous vector field has
order at most 12.
\end{proof}

\begin{corollary}
\label{c1}The vector space $T_{\overline{J_{3}^{v}}(\mathbb{P}^{3}\times
\mathbb{P}^{N_{d}})}\otimes \mathcal{O}_{\mathbb{P}^{3}}(12)\otimes \mathcal{%
O}_{\mathbb{P}^{N_{d}}}(\ast )$ is generated by its global sections on $%
\overline{J_{3}^{v}}(\mathbb{P}^{3}\times \mathbb{P}^{N_{d}})\backslash (\pi
_{3}(\Sigma )\cup \mathcal{X)}.$
\end{corollary}

\begin{remark}
If the third derivative of $f:(\mathbb{C},0)\rightarrow \mathbb{P}^{3}\times
\mathbb{P}^{N_{d}}\backslash \mathcal{X}$ lies inside $\pi _{3}(\Sigma )$
then the image of $f$ is contained in a hyperplane.
\end{remark}

\section{Logarithmic jet differentials}

In this section we recall the basic facts about logarithmic jet
differentials following G. Dethloff and S.Lu \cite{DL96}. Let $X$ be a
complex manifold with a normal crossing divisor $D.$

Let $(X,D)$ be the corresponding complex log-manifold. We start with the
directed manifold $(X,\overline{T}_{X})$ where $\overline{T}_{X}=T_{X}(-\log
D).$ We define $X_{1}:=\mathbb{P(}\overline{T}_{X}),$ $D_{1}=\pi ^{\ast }(D)$
and $V_{1}\subset T_{X_{1}}:$%
\begin{equation*}
V_{1,(x,[v])}:=\{\xi \in \overline{T}_{X_{1},(x,[v])}(-\log D_{1})\text{ };%
\text{ }\pi _{\ast }\xi \in \mathbb{C}v\}
\end{equation*}
where $\pi :X_{1}\rightarrow X$ is the natural projection. If $f:(\mathbb{C}%
,0)\rightarrow (X\backslash D,x)$ is a germ of holomorphic curve then it can
be lifted to $X_{1}\backslash D_{1}$ as $f_{[1]}.$

By induction, we obtain a tower of varieties $(X_{k},D_{k},V_{k})$ with $\pi
_{k}:X_{k}\rightarrow X$ as the natural projection. We have a tautological
line bundle $\mathcal{O}_{X_{k}}(1)$ and we denote $u_{k}:=c_{1}(\mathcal{O}%
_{X_{k}}(1)).$

Let's consider the direct image $\pi _{k\ast }(\mathcal{O}_{X_{k}}(m)).$
It's a locally free sheaf denoted $E_{k,m}\overline{T}_{X}^{\ast }$
generated by all polynomial operators in the derivatives of order $1,2,...,k$
of $f$, together with the extra function $\log s_{j}(f)$ along the $j-th$
component of $D,$ which are moreover invariant under arbitrary changes of
parametrization: a germ of operator $Q\in E_{k,m}\overline{T}_{X}^{\ast }$
is characterized by the condition that, for every germ in $X\backslash D$
and every germ $\phi \in $ $\mathbb{G}_{k}$ of $k$-jet biholomorphisms of $(%
\mathbb{C},0),$%
\begin{equation*}
Q(f\circ \phi )=\phi ^{\prime m}Q(f)\circ \phi .
\end{equation*}

\bigskip

The following theorem makes clear the use of jet differentials in the study
of hyperbolicity:

\bigskip

\textbf{Theorem (\cite{GG80}, \cite{De95}, \cite{DL96}). }\textit{Assume
that there exist integers }$k,m>0$\textit{\ and an ample line bundle }$L$%
\textit{\ on X such that}
\begin{equation*}
H^{0}(X_{k},\mathcal{O}_{X_{k}}(m)\otimes \pi _{k}^{\ast }L^{-1})\simeq
H^{0}(X,E_{k,m}\overline{T}_{X}^{\ast }\otimes L^{-1})
\end{equation*}
\textit{has non zero sections }$\sigma _{1},...,\sigma _{N}.$\textit{\ Let }$%
Z\subset X_{k}$\textit{\ be the base locus of these sections. Then every
entire curve }$f:\mathbb{C}\rightarrow X\backslash D$\textit{\ is such that }%
$f_{[k]}(\mathbb{C})\subset Z.$\textit{\ In other words, for every global }$%
\mathbb{G}_{k}-$\textit{invariant polynomial differential operator P with
values in }$L^{-1},$\textit{\ every entire curve }$f:\mathbb{C}\rightarrow
X\backslash D$\textit{\ must satisfy the algebraic differential equation }$%
P(f)=0.$

\bigskip

If $X\subset \mathbb{P}^{3}$ is a smooth hypersurface$,$ we have
established in \cite{Rou2} the next result:

\bigskip

\label{teo3}\textbf{Theorem (\cite{Rou2}). }\textit{Let X be a
smooth
hypersurface of }$\mathbb{P}^{3}$\textit{\ such that }$d=\deg (X)\geq 92,$%
\textit{\ and A an ample line bundle, then }$E_{3,m}\overline{T}_{\mathbb{P}%
^{3}}^{\ast }\otimes A^{-1}$\textit{\ has global sections for m large enough
and every entire curve }$f:\mathbb{C}\rightarrow \mathbb{P}^{3}\backslash X$%
\textit{\ must satisfy the corresponding algebraic differential equation.}

\bigskip

The proof relies on the filtration of $E_{3,m}\overline{T}_{X}^{\ast }$
obtained in \cite{Rou05}$:$%
\begin{equation*}
Gr^{\bullet }E_{3,m}\overline{T}_{X}^{\ast }=\underset{0\leq \gamma \leq
\frac{m}{5}}{\oplus }(\underset{\{\lambda _{1}+2\lambda _{2}+3\lambda
_{3}=m-\gamma ;\text{ }\lambda _{i}-\lambda _{j}\geq \gamma ,\text{ }i<j\}}{%
\oplus }\Gamma ^{(\lambda _{1},\lambda _{2},\lambda _{3})}\overline{T}%
_{X}^{\ast })
\end{equation*}
where $\Gamma $ is the Schur functor.

This filtration provides a Riemann-Roch computation of the Euler
characteristic \cite{Rou05}:
\begin{eqnarray*}
\chi (\mathbb{P}^{3},E_{3,m}\overline{T}_{\mathbb{P}^{3}}^{\ast })&=&m^{9}(%
\frac{389}{81648000000}d^{3}-\frac{6913}{34020000000}d^{2}\\
&&+\frac{6299}{%
4252500000}d-\frac{1513}{63787500})+O(m^{8}).
\end{eqnarray*}

In dimension 3 there is no Bogomolov vanishing theorem (cf.
\cite{Bo79}) as it is used in dimension 2 to control the
cohomology group $H^{2}$, therefore we need the following
proposition obtained in \cite{Rou2}:

\bigskip

\textbf{Proposition (\cite{Rou2}). }\textit{Let }$\lambda
=(\lambda _{1},\lambda _{2},\lambda _{3})$\textit{\ be a partition
such that }$\lambda _{1}>\lambda _{2}>\lambda _{3}$\textit{\ and
}$\left| \lambda \right| =\sum \lambda _{i}>3d+2.$\textit{\ Then :
}
\begin{equation*}
h^{2}(\mathbb{P}^{3},\Gamma ^{\lambda }\overline{T}_{\mathbb{P}^{3}}^{\ast
})\leq g(\lambda )(d+14)+r(\lambda )
\end{equation*}
\textit{where }$g(\lambda )=\frac{3\left| \lambda \right| ^{3}}{2}\underset{%
\lambda _{i}>\lambda _{j}}{\prod }(\lambda _{i}-\lambda _{j})$\textit{\ and }%
$r$\textit{\ is polynomial in }$\lambda $\textit{\ with homogeneous
components of degrees at most 5.}

\bigskip

This proposition provides the estimate \cite{Rou2}
\begin{equation*}
h^{2}(\mathbb{P}^{3},Gr^{\bullet }E_{3,m}\overline{T}_{\mathbb{P}^{3}}^{\ast
})\leq C(d+14)m^{9}+O(m^{8})
\end{equation*}
where $C$ is a constant.

\section{Proof of theorem \ref{t2}}

Let us consider an entire curve $f:\mathbb{C}\rightarrow \mathbb{P}%
^{3}\backslash X$ for a generic hypersurface of $\mathbb{P}^{3}.$ By
Riemann-Roch and the proposition of the previous section we obtain the
following lemma:

\begin{lemma}
\textit{Let X be a smooth hypersurface of }$\mathbb{P}^{3}$ of degree $d$, $%
0<\delta <\frac{1}{18}$ then $h^{0}(\mathbb{P}^{3},E_{3,m}\overline{T}_{%
\mathbb{P}^{3}}^{\ast }\otimes \overline{K}_{\mathbb{P}^{3}}^{-\delta
m})\geq \alpha (d,\delta )m^{9}+O(m^{8}),$ with
\begin{eqnarray*}
\alpha (d,\delta ) &=&\frac{1}{408240000000}(67737600\delta
^{3}+1945d^{3}-82956d^{2}-968320+ \\
&&1804680d^{2}\delta +12700800d^{2}\delta ^{3}-9408960d^{2}\delta
^{2}+37635840d\delta ^{2}- \\
&&8579520d\delta -50803200d\delta ^{3}-1058400d^{3}\delta
^{3}-105030d^{3}\delta - \\
&&50181120\delta ^{2}+12165120\delta +604704d+784080d^{3}\delta ^{2}).
\end{eqnarray*}
\end{lemma}

\begin{proof}
$E_{3,m}\overline{T}_{\mathbb{P}^{3}}^{\ast }\otimes \overline{K}_{\mathbb{P}%
^{3}}^{-\delta m}$ admits a filtration with graded pieces
\begin{equation*}
\Gamma ^{(\lambda _{1},\lambda _{2},\lambda _{3})}\overline{T}_{\mathbb{P}%
^{3}}^{\ast }\otimes \overline{K}_{\mathbb{P}^{3}}^{-\delta m}=\Gamma
^{(\lambda _{1}-\delta m,\lambda _{2}-\delta m,\lambda _{3}-\delta m)}%
\overline{T}_{\mathbb{P}^{3}}^{\ast }
\end{equation*}
for $\lambda _{1}+2\lambda _{2}+3\lambda _{3}=m-\gamma ;$ $\lambda
_{i}-\lambda _{j}\geq \gamma ,$ $i<j,$ $0\leq \gamma \leq \frac{m}{5}.$

We compute by Riemann-Roch
\begin{equation*}
\chi (\mathbb{P}^{3},E_{3,m}\overline{T}_{\mathbb{P}^{3}}^{\ast }\otimes
\overline{K}_{\mathbb{P}^{3}}^{-\delta m})=\chi (X,Gr^{\bullet }E_{3,m}%
\overline{T}_{\mathbb{P}^{3}}^{\ast }\otimes \overline{K}_{\mathbb{P}%
^{3}}^{-\delta m}).
\end{equation*}
We use the proposition of the previous section to control
\begin{equation*}
h^{2}(X,E_{3,m}\overline{T}_{\mathbb{P}^{3}}^{\ast }\otimes \overline{K}_{%
\mathbb{P}^{3}}^{-\delta m}):
\end{equation*}
\begin{eqnarray*}
h^{2}(\mathbb{P}^{3},\Gamma ^{(\lambda _{1}-\delta m,\lambda _{2}-\delta
m,\lambda _{3}-\delta m)}\overline{T}_{\mathbb{P}^{3}}^{\ast }) &\leq
&g(\lambda _{1}-\delta m,\lambda _{2}-\delta m,\lambda _{3}-\delta m)(d+14)+
\\
&&r(\lambda _{1}-\delta m,\lambda _{2}-\delta m,\lambda _{3}-\delta m)
\end{eqnarray*}
under the hypothesis $\sum \lambda _{i}-3\delta m>3d+2.$ The conditions
verified by $\lambda $ imply $\sum \lambda _{i}\geq \frac{m}{6}$ therefore
the hypothesis will be verified if
\begin{equation*}
m(\frac{1}{6}-3\delta )>3d+2.
\end{equation*}
We conclude with the computation
\begin{equation*}
\chi (\mathbb{P}^{3},E_{3,m}\overline{T}_{\mathbb{P}^{3}}^{\ast }\otimes
\overline{K}_{\mathbb{P}^{3}}^{-\delta m})-h^{2}(\mathbb{P}^{3},Gr^{\bullet
}E_{3,m}\overline{T}_{\mathbb{P}^{3}}^{\ast }\otimes \overline{K}_{\mathbb{P}%
^{3}}^{-\delta m})\leq h^{0}(\mathbb{P}^{3},E_{3,m}\overline{T}_{\mathbb{P}%
^{3}}^{\ast }\otimes \overline{K}_{\mathbb{P}^{3}}^{-\delta m}).
\end{equation*}
\end{proof}

\begin{remark}
If we denote $(\mathbb{P}^{3}\times \mathbb{P}^{N_{d}})_{3}^{v}$ the
quotient of $\overline{J_{3}^{v}}^{reg}(\mathbb{P}^{3}\times \mathbb{P}%
^{N_{d}})$ by the reparametrization group $\mathbb{G}_{3}$, one
can easily verify that each vector field given at section 3
defines a section of the
tangent bundle of the manifold $(\mathbb{P}^{3}\times \mathbb{P}%
^{N_{d}})_{3}^{v}.$
\end{remark}

We have a section
\begin{equation*}
\sigma \in H^{0}(\mathbb{P}^{3},E_{3,m}\overline{T}_{\mathbb{P}^{3}}^{\ast
}\otimes \overline{K}_{\mathbb{P}^{3}}^{-\delta m})\simeq H^{0}((\mathbb{P}%
^{3})_{3},\mathcal{O}_{(\mathbb{P}^{3})_{3}}(m)\otimes \pi _{3}^{\ast }%
\overline{K}_{\mathbb{P}^{3}}^{-\delta m}).
\end{equation*}
with zero set $Z$ and vanishing order $\delta m(d-4).$ Consider the family
\begin{equation*}
\mathcal{X}\subset \mathbb{P}^{3}\times \mathbb{P}^{N_{d}}
\end{equation*}
of hypersurfaces of degree $d$ in $\mathbb{P}^{3}.$ General semicontinuity
arguments concerning the cohomology groups show the existence of a Zariski
open set $U_{d}\subset \mathbb{P}^{N_{d}}$ such that for any $a\in U_{d},$
there exists a divisor
\begin{equation*}
Z_{a}=(P_{a}=0)\subset (\mathbb{P}_{a}^{3})_{3}
\end{equation*}
where
\begin{equation*}
P_{a}\in H^{0}((\mathbb{P}_{a}^{3})_{3},\mathcal{O}_{(\mathbb{P}%
_{a}^{3})_{3}}(m)\otimes \pi _{3}^{\ast }\overline{K}_{(\mathbb{P}%
_{a}^{3})}^{-\delta m})
\end{equation*}
such that the family $(P_{a})_{a\in U_{d}}$ varies holomorphically. We
consider $P$ as a holomorphic function on $\overline{J_{3}}(\mathbb{P}%
_{a}^{3}).$ The vanishing order of this function is no more than $m$ at a
generic point of $\mathbb{P}_{a}^{3}.$ We have $f_{[3]}(\mathbb{C})\subset
Z_{a}.$

Then we invoke corollary \ref{c1} which gives the global generation of
\begin{equation*}
T_{\overline{J_{3}^{v}}(\mathbb{P}^{3}\times \mathbb{P}^{N_{d}})}\otimes
\mathcal{O}_{\mathbb{P}^{3}}(12)\otimes \mathcal{O}_{\mathbb{P}%
^{N_{d}}}(\ast )
\end{equation*}
on $\overline{J_{3}^{v}}(\mathbb{P}^{3}\times
\mathbb{P}^{N_{d}})\backslash (\pi _{3}(\Sigma )\cup
\mathcal{X)}.$

If $f_{[3]}(\mathbb{C)}$ lies in $\pi _{3}(\Sigma )$, $f$ is algebraically
degenerated. So we can suppose it is not the case.

At any point of $f_{[3]}(\mathbb{C)}\backslash \pi _{3}(\Sigma )$ where the
vanishing of $P$ is no more than $m,$ we can find global meromorphic vector
fields $v_{1},...,v_{p}$ $(p\leq m)$ and differentiate $P$ with these vector
fields such that $v_{1}...v_{p}P$ is not zero at this point$.$ From the
above remark, we see that $v_{1}...v_{p}P$ corresponds to an invariant
differential operator and its restriction to $(\mathbb{P}_{a}^{3})_{3}$ can
be seen as a section of the bundle
\begin{equation*}
\mathcal{O}_{(\mathbb{P}_{a}^{3})_{3}}(m)\otimes \mathcal{O}_{\mathbb{P}%
^{3}}(12p-\delta m(d-4)).
\end{equation*}
Assume that the vanishing order of $P$ is larger than the sum of the pole
order of the $v_{i}$ in the fiber direction of $\pi :\mathbb{P}^{3}\times
\mathbb{P}^{N_{d}}\rightarrow \mathbb{P}^{N_{d}}.$ Then the restriction of $%
v_{1}...v_{p}P$ to $\mathbb{P}_{a}^{3}$ defines a jet differential which
vanishes on an ample divisor. Therefore $f_{[3]}(\mathbb{C)}$ should be in
its zero set.

To finish the proof, we just have to see when the vanishing order of $P$ is
larger than the sum of the pole order of the $v_{i}.$ This will be verified
if
\begin{equation*}
\delta (d-4)>12.
\end{equation*}
So we want $\delta >\frac{12}{(d-4)}$ and $\alpha (d,\delta )>0.$ This is
the case for $d\geq 586.$

\noindent \texttt{rousseau@math.u-strasbg.fr}

\noindent D\'{e}partement de Math\'{e}matiques,

\noindent IRMA,\newline
Universit\'{e} Louis Pasteur,

\noindent 7, rue Ren\'{e} Descartes,\newline
\noindent 67084 STRASBOURG CEDEX

\noindent FRANCE

\end{document}